\documentclass[12pt,a4paper,reqno]{amsart}

\usepackage{anysize}
\usepackage[utf8]{inputenc}
\usepackage{amsmath}
\usepackage{amssymb}
\usepackage{array}

\begin{document}

\title[Proofs of some Ramanujan series for $1/\pi$]{Proofs of some Ramanujan series for $1/\pi$ \\ using a Zeilberger's program}

\author{Jes\'{u}s Guillera}

\address{Department of Mathematics, University of Zaragoza, 50009 Zaragoza, SPAIN}

\email{jguillera@gmail.com}

\date{}

\begin{abstract}
We show with some examples how to prove some Ramanujan-type series for $1/\pi$ in an elementary way by using terminating identities. 
\end{abstract}

\maketitle

\subsection*{Introduction}
Up till now, we know how to prove $11$ Ramanujan-type series for $1/\pi$ by using the WZ (Wilf and Zeilberger) method \cite{Gu-on-wz}. Here we will show how to prove some more using a related Zeilberger's algorithm.

\section{The WZ algorithm as a black box}

Let $G(n,k)$ be hypergeometric in $n$ and $k$, that is such that $G(n+1,k)/G(n,k)$ and $G(n,k+1)/G(n,k)$ are rational functions. Then, we can use the Zeilberger's Maple package $\texttt{SumTools[Hypergeometric]);}$.
The output of $\texttt{Zeilberger(G(n,k),k,n,K)[1];}$ is an operator $O(K)$ of the following form
\[
O(K) = P_0(k)+P_1(k) \, K+P_2(k) \, K^2 + \cdots + P_{m}(k) \, K^{m},
\]
where $P_0(k), \, P_1(k), \dots, P_m(k)$ are polynomials of $k$, and $K$ is an operator which increases $k$ in $1$ unity, that is $K G(n,k)= G(n,k+1)$. The output of $\texttt{Zeilberger(G(n,k),k,n,K)[2];}$ gives a function $F(n,k)$ such that
\[
O(K) G(n,k) = F(n+1,k)-F(n,k).
\]
If we sum for $n \geq 0$, we get
\[
O(K) r_k = \lim_{n \to \infty} F(n,k) - F(0,k), \qquad r_k=\sum_{n=0}^{\infty} G(n,k). 
\]
If the above limit and $F(0,k)$ are equal to zero, we have
\[
O(K) r_k = 0,
\]
which is a recurrence of order $m$.

\subsection*{Example 1}
Prove that:
\begin{equation}\label{65-8}
\sum_{n=0}^{\infty} \frac{\left(\frac12\right)_n\left(\frac14\right)_n\left(\frac34\right)_n}{\left(1\right)_n^3}(-1)^n\left(\frac{16}{63}\right)^{2n}(65n+8)=\frac{9\sqrt{7}}{\pi}.
\end{equation}
We have not found a WZ-pair which proves this Ramanujan series. However our proof is closely related to the WZ-method.
\begin{proof}
Let 
\begin{align*}
&A(n,k)=3\left( \frac{64}{63} \right)^k \frac{(-k)_n\left(\frac12\right)_n^2}{\left( \frac12-k\right)_n^2(1)_n}\left( \frac{1}{64}\right)^n(42n+5), 
\\ &
B(n,k)=\frac{(-k)_n\left( \frac{-k}{2} \right)_n\left(\frac12-\frac{k}{2}\right)_n}{\left( \frac12-k\right)_n^2(1)_n}(-1)^n\left(\frac{16}{63}\right)^{2n}(130n-2k+15),
\end{align*}
We define the sequences
\[
r_k=\sum_{n=0}^{\infty} A(n,k), \qquad s_k=\sum_{n=0}^{\infty} B(n,k).
\]
Then we use a Zeilberger's program which finds recurrences. Writing in a Maple session 
\begin{verbatim}
                    with(SumTools[Hypergeometric]);
                    s:=subs(n=0, Zeilberger(A(n,k),k,n,K)[2]);
                    t:=subs(n=0, Zeilberger(B(n,k),k,n,K)[2]);
\end{verbatim}
we see that $s=t=0$. Then, writing
\begin{verbatim}
                    u:=Zeilberger(A(n,k),k,n,K)[1];
                    v:=Zeilberger(B(n,k),k,n,K)[1];
\end{verbatim}
and executing it, we see that $r_k$ and $s_k$ satisfy a common recurrence of order $3$. Then observe that the sums which define $r_k$ and $s_k$ are finite because the terms with $n>k$ are equal to zero due to presence of $(-k)_n$. By direct evaluation, we check that $r_0=s_0$, $r_1=s_1$ and $r_2=s_2$. Hence, as the three first terms are equal, all of them are. Let
\begin{align*}
&r(k)=3\left( \frac{64}{63} \right)^k \sum_{n=0}^{\infty} \frac{(-k)_n\left(\frac12\right)_n^2}{\left( \frac12-k\right)_n^2(1)_n}\left( \frac{1}{64}\right)^n(42n+5), 
\\ &
s(k)=\sum_{n=0}^{\infty} \frac{(-k)_n\left( \frac{-k}{2} \right)_n\left(\frac12-\frac{k}{2}\right)_n}{\left( \frac12-k\right)_n^2(1)_n}(-1)^n\left(\frac{16}{63}\right)^{2n}(130n-2k+15),
\end{align*}
Applying Carlson's theorem \cite[p. 39]{Bailey}, we can deduce that for all complex values of $k$ we have $r(k)=s(k)$. Finally replacing $k=-1/2$, we get
\[
\sum_{n=0}^{\infty} \frac{\left(\frac12\right)_n\left(\frac14\right)_n\left(\frac34\right)_n}{\left(1\right)_n^3}(-1)^n\left(\frac{16}{63}\right)^{2n}(130n+16)=\frac{9\sqrt{7}}{8} \sum_{n=0}^{\infty} \frac{\left(\frac12\right)_n^3}{\left(1\right)_n^3}\left(\frac{1}{64}\right)^n(42n+5)
\]
But in $2002$, we used the WZ-method to prove
\[
\sum_{n=0}^{\infty} \frac{\left(\frac12\right)_n^3}{\left(1\right)_n^3}\left(\frac{1}{64}\right)^n(42n+5)=\frac{16}{\pi},
\]
in an elementary way. Hence we are done.
\end{proof}

\subsection*{Example 2}
Prove that:
\begin{equation}\label{126-10}
\sum_{n=0}^{\infty} \frac{\left(\frac12\right)_n\left(\frac16\right)_n\left(\frac56\right)_n}{\left(1\right)_n^3}\left(\frac{2}{11}\right)^{3n}(126n+10)=\frac{11\sqrt{33}}{2\pi}.
\end{equation}
\begin{proof}
It is completely similar to that in our first example: Use \texttt{Zeilberger} to prove the identity
\begin{multline*}
11  \left(\frac{32}{33}\right)^k \sum_{n=0}^{\infty} \frac{ (-3k)_n\left(\frac13-k\right)_n\left(\frac16-2k\right)_n }{\left( \frac23-2k \right)_n \left( \frac13-4k \right)_n (1)_n } \left( \frac{-1}{8} \right)^n (6n+1)
\\ = \sum_{n=0}^{\infty} \frac{ (-k)_n \left( \frac13-k \right)_n\left(\frac23-k\right)_n}{\left(\frac56-k\right)_n \left(\frac23-2k\right)_n(1)_n }\left( \frac{2}{11} \right)^{3n}(126n+6k+11),
\end{multline*}
and take $k=-1/6$.
\end{proof}

\subsection*{Example 3}
Prove that:
\begin{equation}\label{63-8}
\sum_{n=0}^{\infty} \frac{\left(\frac12\right)_n\left(\frac16\right)_n\left(\frac56\right)_n}{\left(1\right)_n^3}\left(\frac{-4}{5}\right)^{3n}(63n+8)=\frac{5\sqrt{15}}{\pi}.
\end{equation}
\begin{proof}
As in the preceeding examples, first use \texttt{Zeilberger} to show that 
\begin{multline*}
5 \sum_{n=0}^{\infty} \frac{(-3k)_n\left(\frac23+k\right)_n\left(\frac13-k\right)_n }{\left( \frac56-k \right)_n \left( \frac23-2k \right)_n (1)_n} \left( \frac{1}{64} \right)^n (42n+5)
\\ = \left( \frac{15}{16} \right)^{3k} \sum_{n=0}^{\infty} \frac{ (-k)_n \left( \frac13-k \right)_n\left(\frac23-k\right)_n}{\left(\frac56-k\right)_n \left(\frac23-2k\right)_n(1)_n }\left( \frac{-64}{125} \right)^{n}(252n-42k+25).
\end{multline*}
Then take $k=-1/6$.
\end{proof}

\subsection*{Example 4}
With \texttt{Zeilberger}, we can also prove the following general identity:
\[
\sum_{n=0}^{\infty} \frac{(-k)_n\left(\frac12\right)_n^2}{\left( \frac12-k\right)_n^2(1)_n}z^n = (1-z)^k \sum_{n=0}^{\infty} \frac{(-k)_n\left( \frac{-k}{2} \right)_n\left(\frac12-\frac{k}{2}\right)_n}{\left( \frac12-k\right)_n^2(1)_n} \left( \frac{-4z}{(1-z)^2} \right)^n,
\]
which is a particular case of a multi-parameter formula due to Whipple \cite{GS}. Applying to it the operator $5+42 \theta$ at $z=1/64$, where $\theta = z \, d/dz$ (Zudilin's translation method \cite{Zu5}), we get an identity which we have proved directly in Example 1. In a similar way, If we apply the operator $1+6 \theta$ at $z=-1/8$, we get an identity which we can reprove directly with \texttt{Zeilberger}. From this identity we immediatly get an elementary proof of the formula  
\begin{equation}\label{7-1}
\sum_{n=0}^{\infty} \frac{\left(\frac12\right)_n\left(\frac14\right)_n\left(\frac34\right)_n}{\left(1\right)_n^3}\left(\frac{32}{81}\right)^{n}(7n+1)=\frac{9}{2\pi},
\end{equation}
as there is a WZ-method proof of the series in the other side of the identity \cite{Gu-on-wz}.

\subsection*{Example 5}
With \texttt{Zeilberger}, we can also prove the following general identity:
\[
\sum_{n=0}^{\infty}  \frac{(-3k)_n\left(\frac13-k\right)_n\left(\frac16-2k\right)_n }{ \left( \frac23-2k \right)_n \left(\frac13-4k \right)_n (1)_n } z^n 
= 2 \, (4-z)^{3k} \sum_{n=0}^{\infty} \frac{ (-k)_n \left( \frac13-k \right)_n\left(\frac23-k\right)_n}{\left(\frac56-k\right)_n \left(\frac23-2k\right)_n(1)_n}\left( \frac{27z^2}{(4-z)^3} \right)^n,
\]
which is a particular case of a multi-parameter formula due to Bailey \cite{GS}. Applying to it the operator $1+6 \theta$ at $z=-1/8$, we get an identity which we have proved directly in Example 2. In a similar way, if we apply the operators: $1+4 \theta$ at $z=-1$; $1+6 \theta$ at $z=1/4$ and $5+42 \theta$ at $z=1/64$, we get identities which we can reprove directly with \texttt{Zeilberger}. From these identities we can derive respectively the formulas
\begin{align}
\sum_{n=0}^{\infty} \frac{\left(\frac12\right)_n\left(\frac16\right)_n\left(\frac56\right)_n}{\left(1\right)_n^3}\left(\frac35 \right)^{3n}(28n+3) &=\frac{5\sqrt{5}}{\pi}, \label{28-3} \\
\sum_{n=0}^{\infty} \frac{\left(\frac12\right)_n\left(\frac16\right)_n\left(\frac56\right)_n}{\left(1\right)_n^3}\left(\frac{4}{125} \right)^n(11n+1) &=\frac{5\sqrt{15}}{6\pi}, \label{11-1} \\
\sum_{n=0}^{\infty} \frac{\left(\frac12\right)_n\left(\frac16\right)_n\left(\frac56\right)_n}{\left(1\right)_n^3}\left(\frac{4}{85} \right)^{3n}(133n+8) &=\frac{85\sqrt{255}}{54\pi}, \label{133-8}
\end{align}
in an elementary way taking into account that we have shown that they are equal to series that we had already proved by the WZ-method \cite{Gu-on-wz}.

\subsection*{Remarks}
\begin{enumerate}
\item Our proofs are elementary (we do not use the modular theory).
\item Formulas (\ref{11-1}) and (\ref{133-8}) are due to Ramanujan \cite{Ra}. Formulas (\ref{65-8}) and (\ref{7-1}) are due to Berndt, Chan and Liaw \cite{BCL}. Formulas (\ref{126-10}) and (\ref{28-3}) are due to the Borweins \cite{Bo2}. Formula (\ref{63-8}) is due to Baruah and Berndt \cite{BaBe}.
\item For other elementary methods to prove these and other Ramanujan series see \cite{Zu5} and \cite{GuZu}. Those methods are based in the variable $z$, while the proofs in this paper are based in the free parameter $k$.
\end{enumerate}

\end{document}